\documentclass[preprint,3p]{elsarticle}
\usepackage{amsmath}
\usepackage{amssymb}
\newtheorem{thm}{Theorem}[section]
\newtheorem{lem}[thm]{Lemma}

\newtheorem{prop}{Proposition}[section]

\newdefinition{rmk}{Remark}[section]
\newproof{pf}{Proof}
\newproof{pot}{Proof of Theorem \ref{thm2}}
\newproof{poot}{Proof of Corollary \ref{co1}}
\numberwithin{equation}{section}
\newdefinition{ex}{Example}[section]
\journal{}
\begin{document}
\begin{frontmatter}

\title{  Regularity of  extremal solutions of semilinear elliptic problems with non-convex nonlinearities on general domains}
\author{A. AGHAJANI}
\ead{aghajani@iust.ac.ir}.

\address{School of Mathematics, Iran University of Science and Technology, Narmak, Tehran, Iran.}
\address{School of Mathematics, Institute for Research in Fundamental Sciences (IPM), P.O.Box: 19395-5746, Tehran, Iran.}

\begin{abstract}
We consider the semilinear elliptic equation $-\Delta u =\lambda f(u)$
in a smooth bounded domain $\Omega$ of $R^{n}$
with Dirichielt boundary condition, where $f$ is a $C^{1}$ positive and nondeccreasing  function    in $[0,\infty)$ such that $\frac{f(t)}{t}\rightarrow\infty$ as $t\rightarrow\infty$.  When $\Omega$ is an arbitrary domain and $f$ is not necessarily convex, the boundedness of
the extremal solution $u^{*}$ is known only for $n= 2$, established by X. Cabr\'{e} \cite{C1}. In this paper, we prove this  for higher dimensions depending on the nonlinearity $f$. In particular,  we prove that
if
$$\frac{1}{2}<\beta_{-}:=\liminf_{t\rightarrow\infty} \frac{f'(t)F(t)}{f(t)^{2}}\leq \beta_{+}:=\limsup_{t\rightarrow\infty} \frac{f'(t)F(t)}{f(t)^{2}}<\infty$$ where $F(t)=\int_{0}^{t}f(s)ds$, then $u^{*}\in L^{\infty}(\Omega)$, for $n\leq 6$. Also, if $\beta_{-}=\beta_{+}>\frac{1}{2}$ or $\frac{1}{2}<\beta_{-}\leq\beta_{+}<\frac{7}{10}$, then $u^{*}\in L^{\infty}(\Omega)$, for $n\leq 9$. Moreover, if $\beta_{-}>\frac{1}{2}$ then $u^{*}\in H^{1}_{0}(\Omega)$ for  $n\geq 2$.
\end{abstract}

\end{frontmatter}
\section{Introduction}
In this article, we consider the semilinear Dirichlet problem

\begin{equation}
\left\{\begin{array}{ll} -\Delta u= \lambda f(u)& {\rm }\ x\in \Omega,\\~~u>0& {\rm }\ x\in \Omega,\\~~u=0& {\rm }\ x\in \partial \Omega,
\end{array}\right.
\end{equation}
where $\Omega\subset R^{n}$ is a smooth bounded domain, $n\geq1$, $\lambda >0$ is a real parameter, and the nonlinearity $f:[0,\infty]\rightarrow R$ satisfies\\ \\
($H$) ~~~~~~~~~~~$f$ is $C^{1}$, nondecreasing, $f(0)>0$ and $\lim_{s\rightarrow \infty}\frac{f(s)}{s}=\infty.$\\
By a weak solution solution of  $(1.1)$ we mean a nonnegative function $u\in L^{1}(\Omega)$ so that $f(u)\in L_{\delta}^{1}(\Omega)=L^{1}(\Omega,\delta(x)dx)$, $\delta(x)=dist(x,\partial \Omega)$ and
$$\int_{\Omega}(-\Delta\varphi)u=\int_{\Omega}\lambda f(u)\varphi$$
holds for any $\varphi\in C^{2}(\overline{\Omega})$, $\varphi=0$ on $\partial \Omega$ (see Brezis et al. \cite{BCMR}).\\

It is well known (\cite{BCMR,CR,Dup}) that there exists a finite positive extremal parameter $\lambda^{*}$ such that for any $0<\lambda<\lambda^{*}$, problem (1.1) has a minimal classical solution $u_{\lambda}\in C^{2}(\overline{\Omega})$,
while no solution exists, even in the weak sense for $\lambda\geq\lambda^{*}$. The function $\lambda\rightarrow u_{\lambda}$ is increasing and the
increasing pointwise limit $u^{*}(x)=\lim_{\lambda\uparrow\lambda^{*}}u_{\lambda}(x)$ is a weak solution of $(1.1)$ for  $\lambda=\lambda^{*}$
which is called the extremal solution. If $\lambda<\lambda^{*}$ the solution $u_{\lambda}$  is obtained by the implicit function theorem and is stable in the sense that the first Dirichlet
eigenvalue of the linearized problem at $u_{\lambda}$, $-\Delta - \lambda f'(u_{\lambda}$, is positive for all $\lambda\in(0,\lambda_{*})$. That is,
\begin{equation}
\int_{\Omega}|\nabla \varphi|^{2}-\int_{\Omega}\lambda f(u)'\varphi^{2}\geq0,~~~\varphi\in H^{1}_{0}(\Omega).
\end{equation}

The regularity and properties of the extremal solutions have been studied extensively in the literature [2-12,15,19] and it is shown that it depends strongly on the dimension $n$, domain $\Omega$ and nonlinearity $f$. \\

 When $f$ is convex, Nedev in \cite{N} proved that $u^{*}\in L^{\infty}(\Omega)$ for $n=2,3$ in any domain $\Omega$. When $2\leq n \leq 4$  the best
known result was established by Cabr\'{e} \cite{C1} who showed that $u^{*}\in L^{\infty}(\Omega)$ for arbitrary nonlinearity $f$ if in addition $\Omega$ is convex. Applying the  main estimate used
in the proof of the results  of \cite{C1}, Villegas \cite{V} got the same replacing the condition that $\Omega$ is convex with $f$ is convex.
Cabr\'{e} and Capella \cite{CC} proved that $u^{*}\in L^{\infty}(\Omega)$ if $n\leq 9$ and $\Omega=B_{1}$. Also, in \cite{CR}, Cabr\'{e} and Ros-Oton showed that $u^{*}\in L^{\infty}(\Omega)$ if $n\leq 7$ and $\Omega$ is a convex
domain of double revolution (see \cite{CR} for the definition).\\

By imposing extra assumptions on the convex nonlinearity $f$ satisfies (H) much more is known, see \cite{CSS}. Let $f$ is convex and define
\begin{equation}
\tau_{-}:=\liminf_{t\rightarrow\infty} \frac{f(t)f''(t)}{f'(t)^{2}}\leq \tau_{+}:=\limsup_{t\rightarrow\infty} \frac{f(t)f''(t)}{f'(t)^{2}}.
\end{equation}
 Crandall and Rabinowitz \cite{CR} proved $u^{*}\in L^{\infty}(\Omega)$ when $0<\tau_{-}\leq \tau_{+}<2+\tau_{-}+\sqrt{\tau_{-}}$ and $n<4+2\tau_{-}+4\sqrt{\tau_{-}}$. This result was improved by  Ye and Zhou in \cite{YZ} and Sanch\'{o}n in \cite{S1} establishing that $u^{*}\in L^{\infty}(\Omega)$ when $\tau_{-}>0$ and $n<6+4\sqrt{\tau_{-}}$. In \cite{S1} Sanch\'{o}n proved that $u\in L^{\infty}$ whenever $\tau_{-}=\tau_{+}\geq0$ and $n\leq 9$. Recently Cabr\'{e},  Sanch\'{o}n and Spruck \cite{CSS} proved  that if $\tau_{+}<1$ (without assuming $\tau_{-}>0$) and $n<2+\frac{4}{\tau_{+}}$ then $u\in L^{\infty}$, and if $\tau_{+}=1$ and $n<6$ then $u^{*}\in L^{\infty}$. These results improved by the author in \cite{Ag} as follows
 \begin{equation*}
if ~0<\tau_{+}<\infty~~ and~ ~n<\max\{2+\frac{4}{\tau_{+}}+\frac{4}{\sqrt{\tau_{+}}},~4+\frac{2}{\tau_{+}}+\frac{4}{\sqrt{\tau_{+}}}\}~~ then~ ~u^{*}\in L^{\infty}(\Omega).
\end{equation*}
In particular, if $\tau_{+}<\frac{2}{9-2\sqrt{14}}\cong 1.318$ and $n<10$ then $u^{*}\in L^{\infty}(\Omega).$\\

The case when $f$ is not convex and $\Omega$ is arbitrary domain, is more challenging and there is nothing
much in the literature about the boundedness of the extremal solution. Indeed, in this case, again the best result is due to Cabr\'{e} \cite{C1} who showed that $u^{*}\in L^{\infty}(\Omega)$ for arbitrary $f$ and $\Omega$ in dimension $n=2$.\\

In this work we consider problem (1.1) for the case when $f$ is not necessarily convex and $\Omega$ is an arbitrary domain and prove the  boundedness of the extremal solution in higher dimensions under some  extra assumptions on $f$.\\

Let $f$  satisfy (H) and define
\begin{equation}
\beta_{-}:=\liminf_{t\rightarrow\infty} \frac{f'(t)F(t)}{f(t)^{2}}\leq \beta_{+}:=\limsup_{t\rightarrow\infty} \frac{f'(t)F(t)}{f(t)^{2}},
\end{equation}
where $F(t):=\int_{0}^{t}f(s)ds$, for $t\geq0$.\\

The main results of this paper are as follows.
\begin{thm}
Let $f$ (not necessarily convex) satisfy $(H)$ with $\frac{1}{2}<\beta_{-}\leq\beta_{+}<\infty$ and $\Omega$ an arbitrary bounded smooth domain. Let $u^{*}$ be the extremal solution of problem $(1.1)$. Then $u^{*}\in L^{\infty}(\Omega)$ for

\begin{equation}
n<4+4~\Big(\frac{2\beta_{+}-1}{2\beta_{+}}+\sqrt{\frac{2\beta_{-}-1}{\beta_{+}}}\Big).
\end{equation}
Furthermore, if $\beta_{+}<1$ then $u^{*}\in L^{\infty}(\Omega)$ for
\begin{equation}
n<6+\frac{4}{2\beta_{+}-1}\Big(1-\beta_{+}+\sqrt{\beta_{+}(2\beta_{-}-1)}\Big).
\end{equation}
As consequences, by the assumption $\frac{1}{2}<\beta_{-}\leq\beta_{+}<\infty$, we have:\\
(a) If $n\leq 6,$ then $u^{*}\in L^{\infty}(\Omega)$.\\
(b) If $\beta_{-}=\beta_{+}$ or $\beta_{+}<\frac{7}{10}$, then $u^{*}\in L^{\infty}(\Omega)$ for $n\leq 9.$\\
\end{thm}
It is worth mentioning here that, for a convex nonlinearity $f$ we always have $\beta_{+}\geq\beta_{-}\geq \frac{1}{2}$. Indeed in this case $f'$ is a nondecreasing function, hence we have
$$f'(t)F(t)=f'(t)\int_{0}^{t}f(s)ds\geq \int_{0}^{t}f'(s)f(s)ds=\frac{f(t)^{2}}{2}-\frac{f(0)^{2}}{2},$$
now the fact that $f(t)\rightarrow\infty$ as $t\rightarrow\infty$ gives  $\beta_{-}\geq \frac{1}{2}$. Also for general nonlinearities $f$ (not necessarily convex) satisfy only (H) we always have $\beta_{+}\geq \frac{1}{2}$. To see this, by contradiction assume that $0\leq\beta_{+}< \frac{1}{2}$ and take a $\beta\in (\beta_{+},\frac{1}{2})$. Then from the definition of $\beta_{+}$ there exists  $T>0$ such that $f'(t)F(t)\leq\beta f(t)^{2}$ for $t\geq T$, or equivalently,
$\frac{d}{dt}(\frac{f(t)}{F(t)^{\beta}})\leq0$ for $t\geq T$. Thus, $\frac{f(t)}{F(t)^{\beta}}\leq C:=\frac{f(T)}{F(T)^{\beta}}$ for $t\geq T$, and by integration we get $F(t)\leq (C_{1}t+C_{2})^{\frac{1}{1-\beta}}$ for all  $t>T$ and some constants $C_{1},C_{2}$. But, from the superlinearity of $f$ we  have $\lim_{t\rightarrow\infty}\frac{F(t)}{t^{2}}=\infty$, hence we must have $\frac{1}{1-\beta}>2$ or equivalently $\beta>\frac{1}{2}$ which is a contradiction. Hence, we always have $\beta_{+}\geq \frac{1}{2}$.\\
\begin{ex}
Consider problem (1.1) in an arbitrary bounded smooth domain $\Omega$ with  $f(u)=u^{2}+3u+3\cos u+4$. It is easy to see that $f$ satisfies (H), but is not convex (even at infinity). Indeed, we have $f''(u)=2-3\cos u$, which is negative for all $u$ such that $\cos u>\frac{2}{3}$ (so none of the previous results apply). However, by a simple computation we have $1>\beta_{-}=\beta_{+}=\frac{2}{3}>\frac{1}{2}$, hence by Theorem 1.1 we get $u^{*}\in L^{\infty}(\Omega)$ for $n\leq 15.$
\end{ex}
\begin{ex}
As an another example, consider  problem (1.1) with  $f(u)=e^{u}(3+2\cos u)$ and arbitrary bounded smooth domain $\Omega$. Then, $f$ satisfies (H), but is not convex. Indeed, we have $f''(u)=e^{u}(3-4\sin u)$, so $\liminf_{u\rightarrow\infty}f''(u)=-\infty.$ However, we have, after some simplification,
$$\frac{f'(t)(F(t)+4)}{f(t)^{2}}=\frac{(3+2\cos u-\sin u)(3+\sin u+\cos u)}{(3+2\cos u)^{2}}:=\beta(u),$$
which is a periodic function with period $2\pi$, hence (as computed by Mathematica),
$$\beta_{-}=\min_{[0,2\pi]}\beta(u)\approx 0.786244 ~~~ \text{and}~~~ \beta_{+}=\max_{[0,2\pi]}\beta(u)\approx 2.08846,$$
where we used also that $\lim_{t\rightarrow\infty}\frac{4f'(t)}{f(t)^{2}}=0$. Now, using Theorem (1.1) we get $u^{*}\in L^{\infty}(\Omega)$ for $n\leq 9.$
\end{ex}

Now consider the well-known convex nonlinearities $f(t)=e^{t}$ or $(1+t)^{p}$, $p>1$. When $f(t)=e^{t}$ we have $\beta_{-}=\beta_{+}=1$ then from Theorem 1.1 we get $u^{*}\in L^{\infty}(\Omega)$ for $n\leq 9.$ Also, for $f(t)=(1+t)^{p}$, $p>1$ we have $\beta_{-}=\beta_{+}=\frac{p}{p+1}<1$, hence from Theorem 1.1 we get
$$u\in L^{\infty}(\Omega)~~for~~n<2\Big(1+\frac{2p}{p-1}+2\sqrt{\frac{p}{p-1}}\Big).$$
The above results are well-known in the literature \cite{CR, YZ, S1}.\\

Also, notice that if $f$ is convex and $\tau_{-}>0$ then we must have $\beta_{-}>\frac{1}{2}$. Indeed, it is easy to see that the condition $\tau_{-}>0$ yields that for every $0<\tau<\tau_{-}$, there exists  $T=T(\tau)>0$ such that the function $\frac{f'(t)}{f(t)^{\tau}}$ is increasing in $[T,\infty)$, hence
$$f'(t)F(t)=f(t)^{\tau}\frac{f'(t)}{f(t)^{\tau}}\int_{0}^{t}f(s)ds\geq f(t)^{\tau} \int_{T}^{t}\frac{f'(s)}{f(s)^{\tau}}f(s)ds=\frac{f(t)^{2-\tau}}{2-\tau}-\frac{f(T)^{2-\tau}}{2-\tau},$$
now the facts that $f(t)\rightarrow\infty$ as $t\rightarrow\infty$ and $\tau>\tau_{-}$ was arbitrary give $\beta_{-}\geq \frac{1}{2-\tau_{-}}>\frac{1}{2}$.\\

To get the regularity of the extremal solution in low dimensions or proving that it is in the energy class (i.e., $u^{*}\in H^{1}_{0}(\Omega)$) we can  weaken the assumptions as follows.

\begin{thm}
Let $f$ (not necessarily convex) satisfy $(H)$ and $\Omega$ an arbitrary bounded smooth domain in $R^{n}$. Let $u^{*}$ be the extremal solution of problem $(1.1)$. Then\\
(i) if for some $\epsilon>0$ there exist $t_{0}>0$ such that we have
\begin{equation}
\frac{f'(t)F(t)}{f(t)^{2}}\geq\frac{1}{2}+\frac{\epsilon t}{f(t)},~~t>t_{0}
\end{equation}
then $u^{*}\in H^{1}_{0}(\Omega)$. In particular this is true if $\beta_{-}>\frac{1}{2}$.\\
(ii) If for some $0<\delta\leq1$, $\lim_{t\rightarrow\infty}\frac{f(t)}{t^{2-\delta}}=\infty$, and there exist $t_{0}>0$ such that
\begin{equation}
\frac{f'(t)F(t)}{f(t)^{2}}\geq\frac{1}{2}+\frac{1 }{t^{2-\delta}},~~t>t_{0}
\end{equation}
then $u^{*}\in L^{\infty}(\Omega)$ for $n<5$. In particular this is true if $\beta_{-}>\frac{1}{2}$.
\end{thm}
\section{Preliminaries and Auxiliary Results}

To prove the main results we need the following simple technical lemma based on inequality (1.2), which is used frequently in the literature, for example \cite{CR,CSS,N,YZ}. It is also proved in \cite{Ag} for the general semilinear elliptic equation $-Lu=\lambda f(u)$ with zero Dirichlet boundary condition, but for the convenience of
the reader we sketch a proof here for the case $L=\Delta$. Also,
\begin{lem}
Let $u_{\lambda}$  be the minimal solution of $(1.1)$ and $g:[0,\infty]\rightarrow [0,\infty]$ be a $C^{1}$ function with $g(0)=0$ and satisfy
\begin{equation}
H(t):=g(t)^{2}f'(t)-G(t)f(t)\geq0,~~\text{for}~t~\text{sufficiently~large},
\end{equation}
where $G(t):=\int_{0}^{t}g'(s)^{2}ds$. Then $||H(u_{\lambda})||_{ L^{1}(\Omega)}\leq C$, where $C$ is a constant independent of $\lambda$.
\end{lem}
\begin{pf}
Let $u_{\lambda}\in C^{2}(\overline{\Omega})$  be the minimal classical solution of (1.1) where $0<\lambda<\lambda^{*}$, and take $\varphi=g(u_{\lambda})$ in the semi-stability condition (1.2). Then we get
\begin{equation}
\int_{\Omega}g'(u_{\lambda})^{2}|\nabla u_{\lambda}|^{2}dx-\int_{\Omega}\lambda f'(u_{\lambda})g(u_{\lambda})^{2}dx\geq0.
\end{equation}
By using the Green's formula one can show that
\begin{equation}
\int_{\Omega}g'(u_{\lambda})^{2}|\nabla u_{\lambda}|^{2}dx=\int_{\Omega}\lambda G(u_{\lambda})f( u_{\lambda})dx.
\end{equation}
Using (2.3) in (2.2) we obtain
\begin{equation}
\int_{\Omega} H(u_{\lambda})dx\leq0.
\end{equation}
Now from (2.1) there exists $t_{0}>0$ so that $H(t)\geq0$ for $t\geq t_{0}$, thus using (2.4) we obtain
$$\int_{\Omega} |H(u_{\lambda})|dx=\int_{u_{\lambda}\leq t_{0}} |H(u_{\lambda})|dx+\int_{u_{\lambda}\geq t_{0}} H(u_{\lambda})dx$$
$$\leq \int_{u_{\lambda}\leq t_{0}} (|H(u_{\lambda})|-H(u_{\lambda}))dx\leq C_{0}|\Omega|,$$
where $|\Omega|$ denotes the Lebesgue measure of $\Omega$ and $C_{0}:=\sup_{t\in[0,t_{0}]}(|H(t)|-H(t))$. Now, since $C_{0}$ is independent of $\lambda$ we get the desired result. $\blacksquare$
\end{pf}
The following consequence of the above lemma is essential in the proof of the main results.
\begin{prop}
Let $u_{\lambda}$  be the minimal solution of $(1.1)$ and $\xi:[0,\infty]\rightarrow [0,\infty]$ be a $C^{1}$ function such that for some $t_{0}>0$ we have $\xi(t)\leq \frac{f'(t)}{f(t)}$, $\xi'(t)+\xi(t)^{2}\geq0$ for $t\geq t_{0}$, and $\frac{E(t)}{f(t)}\rightarrow\infty$ as $t\rightarrow\infty$ where
\begin{equation}
E(t):=f(t)~\Big(\frac{f'(t)}{f(t)}-\xi(t)\Big)~e^{2\int_{t_{0}}^{t}(\xi(s)+\sqrt{\xi'(s)+\xi(s)^{2}})ds}.
\end{equation}
Then $||E(u_{\lambda})||_{ L^{1}(\Omega)}\leq C$, where $C$ is a constant independent of $\lambda$.
\end{prop}
\begin{pf}
Let $g:[0,\infty]\rightarrow [0,\infty]$ be a $C^{1}$ function with $g(0)=0$ and
$$g(t)=e^{\int_{t_{0}}^{t}(\xi(s)+\sqrt{\xi'(s)+\xi(s)^{2}})ds},~~ \text{for} ~t\geq t_{0}.$$
Also, let $G(t)=\int_{0}^{t}g'(s)^{2}ds$ as in lemma 2.1. Then using the equality
$$g'(t)=(\xi(t)+\sqrt{\xi'(t)+\xi(t)^{2}})g(t)~~ ~\text{for}~~~t\geq t_{0},$$
we compute
$$\frac{d}{dt}\Big(\xi(t)g(t)^{2}-G(t)\Big)=\xi'(t)g(t)^{2}+2\xi(t)g(t)g'(t)-g'(t)^{2}$$
$$=g(t)^{2}\Big(\xi'(t)+2\xi(t)^{2}+2\xi(t)\sqrt{\xi'(t)+\xi(t)^{2}}\Big)-g'(t)^{2}$$
$$=g(t)^{2}\Big(\xi(s)+\sqrt{\xi'(t)+\xi(t)^{2}}\Big)^{2}-g'(t)^{2}$$
$$=g'(t)^{2}-g'(t)^{2}=0,~~ \text{for} ~t\geq t_{0},$$
implies that
\begin{equation}
G(t)=\xi(t)g(t)^{2}+C_{0},~\text{where}~C_{0}:=G(t_{0})-\xi(t_{0}).
\end{equation}
Now using (2.6), for $t\geq t_{0}$ we have
$$H(t):=g(t)^{2}f'(t)-G(t)f(t)=g(t)^{2}f'(t)-g(t){2}\xi(t)f(t)-C_{0}f(t)=E(t)-C_{0}f(t),$$
which is positive for large $t\geq t_{0}$ (by the assumption), hence by Lemma 2.1 we get $||H(u_{\lambda})||_{ L^{1}(\Omega)}\leq C_{1}$, where $C_{1}$ is a constant independent of $\lambda$. However, again by the assumption we have $0<E(t)<2H(t)$ for large $t$ that also gives $||E(u_{\lambda})||_{ L^{1}(\Omega)}\leq C_{2}$, where $C_{2}$ is a constant independent of $\lambda$, which is the desired result. $\blacksquare$
\end{pf}
To prove Theorem 1.2 in the next section, we also need the following rather standard result. For a simple proof see \cite{Ag}.
\begin{prop}
Let $f$ satisfy (H) and $u_{\lambda}$ be the minimal solution of problem (1.1). If there exists a positive constant $C$ independent of $\lambda$ such that
\begin{equation}
||u_{\lambda}||_{L^{1}(\Omega)}\leq C ~~and~~||\frac{\tilde{f}(u_{\lambda})^{\alpha}}{u_{\lambda}^{\sigma}}||_{L^{1}(\Omega)}\leq C,~~for~some~0\leq \sigma\leq\alpha,
\end{equation}
where $\tilde{f}(u)=f(u)-f(0)$ and $\alpha\geq 1$, then
\begin{equation}
||u_{\lambda}||_{L^{\infty}(\Omega)}\leq \tilde{C}~~for~~n<2\alpha,
\end{equation}
where $\tilde{C}$ is a positive constant independent of $\lambda$.
\end{prop}
\section{Proof of the main results}
{\bf Proof of Theorem 1.1}\\
By the assumptions we have $\frac{1}{2}<\beta_{-}\leq\beta_{+}<\infty$. Take  $\frac{1}{2}<\beta_{1}<\beta_{2}<\beta_{-}$ and  $\beta_{3}\in(\beta_{+},\infty)$, then by the definition of $\beta_{-},\beta_{+}$ (see (1.4)) there exists a $t_{0}>0$ such that
\begin{equation}
 \beta_{1}<\beta_{2}< \frac{f'(t)F(t)}{f(t)^{2}}<\beta_{3},~~\text{for}~~t\geq t_{0}.
\end{equation}
Now let $\xi:[0,\infty]\rightarrow [0,\infty]$ be a $C^{1}$ function such that  $\xi(t)=\beta_{1}\frac{f(t)}{F(t)}$ for $t\geq t_{0}$, where $F(t):=\int_{0}^{t}f(s)ds$. Then from (3.1) we have $\xi(t)\leq \frac{f'(t)}{f(t)}$ and
\begin{equation}
\xi'(t)+\xi(t)^{2}=\beta_{1}\Big(\frac{f'(t)}{F(t)}-(1-\beta_{1})\frac{f(t)^{2}}{F(t)^{2}}\Big)\geq (2\beta_{1}-1)\frac{f'(t)}{F(t)}
\geq \frac{2\beta_{1}-1}{\beta_{3}}\frac{f'(t)^{2}}{f(t)^{2}},
\end{equation}
for $t\geq t_{0}$. Notice that, in (3.2) we used the fact that $\beta_{1}<1$ (because we always have $\beta_{-}\leq1$). Also, from (3.1) we have $$\frac{f'(t)}{f(t)}-\xi(t)\geq (\beta_{2}-\beta_{1})\frac{f(t)}{F(t)}, ~~\text{for}~~t\geq t_{0}.$$
Now let the function $E(t)$ be given as in (2.5) in Proposition 2.1. By the later inequality,  (3.1), (3.2) and the fact that
$\int_{t_{0}}^{t}\xi(s)ds=\beta_{1}(\ln F(t)-\ln F(t_{0}))$, we have
\begin{equation}
E(t)=f(t)\Big(\frac{f'(t)}{f(t)}-\xi(t)\Big)e^{2\int_{t_{0}}^{t}(\xi(s)+\sqrt{\xi'(s)+\xi(s)^{2}})ds}\geq CF(t)^{2\beta_{1}-1}
f(t)^{2+2\sqrt{\frac{2\beta_{1}-1}{\beta_{3}}}},
\end{equation}
where $ C$ is a positive constant depends only on $f$. Now, writing the last inequality in (3.1) as $\frac{f'(t)}{f(t)}<\beta_{3}\frac{f(t)}{F(t)}$
for $t_{0}>0$, then integration  from $t_{0}$ to $t$ gives
 \begin{equation}
F(t)\geq C f(t)^{\frac{1}{\beta_{3}}}~~\text{ for}~~t\geq t_{0}.
\end{equation}
Using (3.4) in (3.3) we arrive at
\begin{equation}
E(t)\geq
f(t)^{\gamma},~~where ~~\gamma:=2+\frac{2\beta_{1}-1}{\beta_{3}}+2\sqrt{\frac{2\beta_{1}-1}{\beta_{3}}},~\text{for} ~t\geq t_{0}.
\end{equation}
And, since $\frac{E(t)}{f(t)}\rightarrow\infty$ $t\rightarrow\infty$, from Proposition 2.1 we get $||E(u_{\lambda})||_{ L^{1}(\Omega)}\leq C$ and then from (3.5) $||f(u_{\lambda})^{\gamma}||_{ L^{1}(\Omega)}\leq C$, where $C$ is a constant independent of $\lambda$. Now the standard elliptic regularity theory gives $u^{*}\in L^{\infty}(\Omega)$ for $n<2\gamma$, and since $\beta_{1}$ and $\beta_{3}$ were arbitrary in the intervals $(\frac{1}{2},\beta_{-})$ and $(\beta_{+},\infty)$, respectively, thus $u^{*}\in L^{\infty}(\Omega)$ for
\begin{equation}
n<4+4~\Big(\frac{2\beta_{+}-1}{2\beta_{+}}+\sqrt{\frac{2\beta_{-}-1}{\beta_{+}}}\Big):=\gamma_{1},
\end{equation}
that proves the first part.\\
Now assume that $\beta_{+}<1$, then we can also assume that $\beta_{3}<1$. Now, from (3.4) we have $f(t)F(t)^{-\beta_{3}}\leq C_{1}$, for $t\geq t_{0}$, and  integration from $t_{0}$ to $t$ gives $F(t)\leq C_{2}t^{\frac{1}{1-\beta_{3}}}$, $t\geq t_{1}$, for some $t_{1}\geq t_{0}$. This together (3.4) implies that $f(t)\leq C_{3}t^{\frac{\beta_{3}}{1-\beta_{3}}}$, for $t\geq t_{1}$, that also  yields, for some  $t_{2}\geq t_{1}$,
$$f(t)^{\gamma_{1}}\geq C_{4}\frac{f(t)^{\gamma_{2}}}{t^{\gamma_{2}}},~~\text{ for}~~t\geq t_{2},~~\gamma_{2}:=\frac{\beta_{3}}{2\beta_{3}-1}\gamma_{1},$$
where $\gamma_{1}$ is given in (3.6). Hence,  $||\frac{\tilde{f}(u_{\lambda})^{\gamma_{2}}}{u_{\lambda}^{\gamma_{2}}}||_{ L^{1}(\Omega)}\leq C$, where $\tilde{f}(t)=f(t)-f(0)$ and $C$ is a constant independent of $\lambda$. Now, from Proposition 2.2 gives
$u^{*}\in L^{\infty}(\Omega)$ for $n<2\gamma_{2}$, that proves the second part.\\
To prove part (a), note that in the case $\beta_{+}\geq 1$, it is easy to see that the right hand side of (1.5) is larger than $6$ and  when  $\beta_{+}<1$ we can use (1.6).\\
To show part (b), note that (from (1.6)) if $6+\frac{4(1-\beta_{+})}{2\beta_{+}-1}>9$, which is equivalent to $\beta_{+}<\frac{7}{10}$, then $u^{*}\in L^{\infty}(\Omega)$ for $n<9$.\\
Now assume that $\beta_{-}=\beta_{+}$ then we have $\beta_{+}\leq 1$ (as we always have $\beta_{-}\leq1$). Also, from the later part we can consider only the case $\beta_{+}\geq\frac{7}{10}$.  If  $\beta_{-}=\beta_{+}=1$ then from (1.5) we have $u^{*}\in L^{\infty}(\Omega)$ for $n<10$. Also, if $\frac{7}{10}\leq\beta_{+}<1$ we  can use (1.6). We need to show that the right hand side of (1.6) is larger than $9$. In the case $\frac{7}{10}\leq\beta_{+}<1$ this is equivalent to $68\beta_{+}^{2}+49<124\beta_{+}$, which  obviously holds for $\beta_{+}\in[\frac{7}{10},1)$. $\blacksquare$\\
{\bf Proof of Theorem 1.2}\\
Let $\xi:[0,\infty]\rightarrow [0,\infty]$ be a $C^{1}$ function such that  $\xi(t)=\frac{1}{2}\frac{f(t)}{F(t)}$ for $t\geq t_{0}$, where $t_{0}$ given as in (1.7). Then from (1.7) we have
\begin{equation}
\frac{f'(t)}{f(t)}-\xi(t)=\frac{f(t)}{F(t)}\Big(\frac{f'(t)F(t)}{f(t)^{2}}-\frac{1}{2}\Big)\geq\epsilon \frac{t}{F(t)},~~ \text{for}~~t\geq t_{0},
\end{equation}
and
\begin{equation}
\xi'(t)+\xi(t)^{2}=\frac{1}{2}\Big(\frac{f'(t)F(t)}{f(t)^{2}}-\frac{1}{2}\Big)\frac{f(t)^{2}}{F(t)^{2}}>0~\text{for} ~t\geq t_{0}.
\end{equation}
From (3.7) and (3.8) we get
\begin{equation}
E(t)=f(t)\Big(\frac{f'(t)}{f(t)}-\xi(t)\Big)e^{2\int_{t_{0}}^{t}(\xi(s)+\sqrt{\xi'(s)+\xi(s)^{2}})ds}\geq Ctf(t),~\text{for} ~t\geq t_{0}
\end{equation}
where $ C$ is a positive constant depends only on $f$. Thus, from Proposition 2.1 we get $u_{\lambda}f(u_{\lambda})\in L^{1}(\Omega)$ by a constant independent of $\lambda$. Multiplying (1.1) by $u_{\lambda}$ we get
$$\int_{\Omega}|\nabla u_{\lambda}|^{2}dx=\lambda\int_{\Omega}u_{\lambda}f(u_{\lambda})dx\leq\lambda^{*}C,$$
with $C$ independent of $\lambda$, which leads to $\int_{\Omega}|\nabla u^{*}|^{2}dx\leq \lambda^{*}C$.
Hence, $u^{*}\in H^{1}_{0}(\Omega)$, that proves part (i).\\
Let $\xi(t)$ be as defined above. Then, using (1.8) and similar to the proof of part (i) we can show that $E(t)\geq C\frac{f(t)^{2}}{t^{2-\delta}}$ for $t\geq t_{0}$, hence by the assumption  $\frac{E(t)}{f(t)}\rightarrow\infty$ as $t\rightarrow\infty$. Now from Proposition 2.1 we get
\begin{equation}
\frac{\tilde{f}(u_{\lambda})^{2}}{u_{\lambda}^{2-\delta}}\in L^{\infty}(\Omega).
\end{equation}
 Now we proceed  similar to the proof of Theorem 1.1 in \cite{V}. From Proposition 2.3 in \cite{V}, there exits a universal constant $C_{1}$ independent of $f$ , $\Omega$  and $\lambda$  such that
  \begin{equation}
||u_{\lambda}||_{ L^{\infty}(\Omega)}\leq C_{1}||\nabla u_{\lambda}||_{ L^{4}(\Omega)}.
\end{equation}
 Also, from the continuous inclusion $W^{2,2}(\Omega)\subset W^{1,4}(\Omega)$ and elliptic regularity theory (see \cite{ADN}), there exits  a constant $C_{2}=C_{2}(\Omega)$   such that
 \begin{equation}
||u_{\lambda}||_{ W^{1,4}(\Omega)}\leq C_{2}||u_{\lambda}||_{ W^{2,2}(\Omega)}~~~\text{and}~~~||u_{\lambda}||_{ W^{2,2}(\Omega)}\leq C_{2}||\lambda f(u_{\lambda})||_{ L^{2}(\Omega)}.
\end{equation}
Now from (3.10), (3.11) and (3.12)  we have (in the following inequalities various constants will be denoted by $C$)
$$||u_{\lambda}||_{ L^{\infty}(\Omega)}\leq C||\nabla u_{\lambda}||_{ L^{4}(\Omega)}\leq C|| u_{\lambda}||_{ W^{1,4}(\Omega)}\leq C|| u_{\lambda}||_{ W^{2,2}(\Omega)}$$
$$\leq C||\lambda f(u_{\lambda})||_{ L^{2}(\Omega)}\leq C\lambda^{*}\Big(\int_{u_{\lambda}\leq1}f(u_{\lambda})^{2}dx+\int_{u_{\lambda}>1}f(u_{\lambda})^{2}dx\Big)^{\frac{1}{2}}$$
$$\leq C\Big(f(1)^{2}|\Omega|+\int_{u_{\lambda}>1}\frac{f(u_{\lambda})^{2}}{u_{\lambda}^{2-\delta}}u_{\lambda}^{2-\delta}dx\Big)
^{\frac{1}{2}}\leq  C\Big(f(1)^{2}|\Omega|+\tilde{C}||u_{\lambda}||^{2-\delta}_{L^{\infty}(\Omega)}\Big)
^{\frac{1}{2}},$$
for every $\lambda\in(0,\lambda_{*})$. Hence, we must have
$$||u_{\lambda}||^{2}_{ L^{\infty}(\Omega)}\leq A+B||u_{\lambda}||^{2-\delta}_{L^{\infty}(\Omega)},~\text{for ~every} \lambda\in(0,\lambda_{*}),$$
where $A$ and $B$ are positive constants independent of $\lambda$. This implies that $||u_{\lambda}||_{ L^{\infty}(\Omega)}\leq C$ with $C$ independent of $\lambda,$ now letting $\lambda\rightarrow \lambda^{*}$ gives $u^{*}\in L^{\infty}(\Omega)$. $\blacksquare$
\section{Acknowledgement}
This research was in part supported by a grant from IPM (No. 93340123).

\end{document}